\documentclass[11pt,a4paper]{article}
\usepackage[final]{optional}

\usepackage{amsfonts, amsmath, wasysym}
\usepackage{stmaryrd} 

\usepackage{tikz}
\usetikzlibrary{shapes, positioning}

\usepackage{amssymb,amsthm,
paralist
}

\usepackage{
latexsym,
}

\usepackage{url}

\definecolor{darkgreen}{rgb}{0,0.5,0}
\definecolor{darkred}{rgb}{0.7,0,0}
\usepackage[colorlinks, 
citecolor=darkgreen, linkcolor=darkred
]{hyperref}

\theoremstyle{plain}

\numberwithin{equation}{section}

\newcommand{\s}{\ensuremath{{\mathcal S}}}

\newcommand{\ci}{\ensuremath{{\mathcal I}}}

\newcommand{\calr}{\ensuremath{{\mathcal R}}}

\newcommand{\pl}[2]{{\frac{\partial #1}{\partial #2}}}

\newcommand{\al}{\alpha}

\newcommand{\de}{\delta}
\newcommand{\om}{\omega}
\newcommand{\Om}{\Omega}

\newcommand{\la}{\lambda}

\newcommand{\Si}{\Sigma}

\newcommand{\ep}{\varepsilon}

\newcommand{\R}{\ensuremath{{\mathbb R}}}
\newcommand{\N}{\ensuremath{{\mathbb N}}}

\newcommand{\C}{\ensuremath{{\mathbb C}}}

\newcommand{\downto}{\downarrow}
\newcommand{\upto}{\uparrow}

\newcommand{\lap}{\Delta}

\DeclareMathOperator{\VolB}{VolB}

\DeclareMathOperator{\inj}{inj}

\newcommand{\beq}{\begin{equation}}
\newcommand{\eeq}{\end{equation}}
\newcommand{\beqa}{\begin{equation}\begin{aligned}}
\newcommand{\eeqa}{\end{aligned}\end{equation}}
\newcommand{\brmk}{\begin{rmk}}
\newcommand{\ermk}{\end{rmk}}
\newcommand{\partref}[1]{\hbox{(\csname @roman\endcsname{\ref{#1}})}}

\newcommand{\Ric}{{\mathrm{Ric}}}
\newcommand{\Scal}{{\mathrm{Scal}}}

 \newtheorem{thm}{Theorem}[section]

\newtheorem{lem}[thm]{Lemma}

\newtheorem{defn}[thm]{Definition}
\newtheorem{rmk}[thm]{Remark}
\newtheorem{conj}[thm]{Conjecture}

\newcommand{\CPICo}{{\ensuremath{\mathrm{C_{PIC1}}}}}
\newcommand{\CPICt}{\ensuremath{\mathrm{C_{PIC2}}}}
\newcommand{\cone}{\ensuremath{\mathrm{C}}}
\newcommand{\CB}{\ensuremath{\mathcal{C}_B(\R^n)}}

\title{\sc RICCI FLOW AND PIC1}
\author{Peter M. Topping}
\date{1 September 2023}

\begin{document}

%

\maketitle

\begin{abstract}
We survey several problems concerning Riemannian manifolds with positive curvature of one form or another. We describe the PIC1 notion of positive curvature and argue that it is often the sharp notion of positive curvature to consider. Finally we explain how recent Ricci flow theory is particularly well adapted to solve these problems.
\end{abstract}

\tableofcontents

\section {Preamble}

One of the great successes of Ricci flow since its inception is to enable us to draw global conclusions about a manifold from local geometric hypotheses. In particular, given a Riemannian manifold satisfying a pointwise constraint on its curvature, Ricci flow can often be used to collate this pointwise information to give global topological information such as the diffeomorphism type of the manifold.
In this survey we take a look at several instances of this principle and 
make the case that the most natural curvature hypotheses in many of these situations involves the so-called PIC1 curvature condition that was introduced by
Micallef and Moore in 1988 \cite[Section 5]{MM}.
We give an introduction to the curvature condition from scratch and survey a collection of results and conjectures that concern it. In addition to the local-global class of results just alluded to, we will see that the PIC1 condition is a particularly natural condition to impose on a sequence of manifolds in order to control the Gromov-Hausdorff limit.

\section{Notions of positive curvature}

We begin by surveying the main notions of positive (or nonnegative) curvature that will concern us.
These notions, and their relationships, are best described at the level of algebraic curvature tensors. Working in $\R^n$, for $n\geq 2$, 
we denote by 
$\mathcal{C}_B(\R^n)$ the vector space of symmetric bilinear forms on 
$\Lambda^2\R^n$ that satisfy the Bianchi identity. 
We adopt the convention that given $\calr\in \mathcal{C}_B(\R^n)$ and orthonormal $e_1,e_2\in \R^n$, the sectional curvature of the plane spanned by $e_1$ and $e_2$ 
is $\calr(e_1\wedge e_2,e_1\wedge e_2)$.

The Euclidean metric $g$ on $\R^n$ extends to tensor products of $\R^n$, and in particular to $\Lambda^2\R^n$. Up to normalisation convention, this extended inner product coincides with the curvature tensor $\ci\in \mathcal{C}_B(\R^n)$ of constant sectional curvature $1$. Using $\ci$ as a metric, we view elements of $\mathcal{C}_B(\R^n)$ equivalently as endomorphisms of $\Lambda^2\R^n$, with $\ci$ itself being the identity.

We will consider nonnegative curvature conditions generally as closed, convex, $O(n)$-invariant cones within $\mathcal{C}_B(\R^n)$. That way, when confronted with a Riemannian manifold, we can say it satisfies the curvature condition if the restriction of the curvature tensor to each tangent space, identified with $\R^n$ using the metric, satisfies the same condition. The $O(n)$-invariance ensures that the notion is independent of the identification. For each of our cones, the corresponding \emph{positive} curvature condition will be that the curvature tensor should lie in the \emph{interior} of the cone.

The simplest and strongest nonnegative curvature condition that we will consider here
is the familiar notion of nonnegative curvature operator. We say that $\calr\in \mathcal{C}_B(\R^n)$
satisfies this condition if it is nonnegative definite as a bilinear form on
$\Lambda^2\R^n$, i.e. all the eigenvalues of $\calr$ are nonnegative. 

We will be mainly interested in weaker curvature conditions. A little weaker is the 
notion of 2-nonnegative curvature operator. Traditionally $\calr$ is said to satisfy this condition if the sum of the two lowest eigenvalues of $\calr$ is nonnegative, but we will adopt a different viewpoint shortly that emphasises why it is particularly significant to add two eigenvalues rather than any larger number of eigenvalues.

To go further, we must complexify. Each $\calr\in \mathcal{C}_B(\R^n)$ can be extended by complex linearity to a symmetric bilinear form on $\Lambda^2\C^n$.
An equivalent definition of \textbf{nonnegative curvature operator} is then that 
$$\calr(\om,\overline{\om})\geq 0\quad\text{ for all }\om\in \Lambda^2\C^n.$$
A key notion is that of complex sectional curvature, which appears 
for Riemannian manifolds in the work of Micallef-Moore \cite{MM}.
We say that $\calr\in \mathcal{C}_B(\R^n)$ has \textbf{nonnegative complex sectional curvature} (also known as \textbf{WPIC2}, for reasons that we will mention later) if 
$$\calr(\om,\overline{\om})\geq 0\quad\text{ for all \textbf{simple} }\om\in \Lambda^2\C^n.$$
Here $\om$ is said to be simple if there exist $v,w\in\C^n$ such that $\om=v\wedge w\neq 0$, i.e. if it is of rank one.
We refer to the complex subspace $\Si$ of $\C^n$ spanned by $v,w$ as a complex section.
We can define the complex sectional curvature $K^\C(\Si)\in\R$ of $\Si$ by 
$\calr(\om,\overline{\om})$ if we choose $v,w$ above to be orthonormal with respect to the Hermitian inner product on $\C^n$, but it is generally wise to discuss only notions that involve the complex linear extension of $g$ to $\C^n$ (also denoted $g$) and avoid 
notions that are explicitly reliant on the Hermitian inner product defined by 
$\langle v,w\rangle:=g(v,\overline{w})$, for reasons that will become clear.
We will denote the closed convex $O(n)$-invariant cone of all $\calr\in\CB$
satisfying the WPIC2 condition as $\CPICt$. A manifold is said to be PIC2 if its 
curvature tensor lies in the interior of $\CPICt$, or equivalently that it has positive complex sectional curvature.

For many purposes, the condition of nonnegative complex sectional curvature is  stronger than necessary. One effective way to weaken it is to require only nonnegativity of the 
complex sectional curvature corresponding to a restricted subset of all possible complex sections. For example, we could restrict to all real sections, by which we mean all $\Si$ such that $\Si=\overline \Si$. That would simply recover the classical notion of nonnegative sectional curvature, which is not so significant for the present survey. 
A more interesting restricted class arises by considering only \textit{totally isotropic} sections $\Si$, i.e. sections for which every $v\in \Si$ is isotropic.
Recall that $v\in \C^n$ is said to be isotropic if $g(v,v)=0$.
We say that $\calr\in \mathcal{C}_B(\R^n)$ has \textbf{nonnegative isotropic curvature}
(or weakly positive isotropic curvature, abbreviated \textbf{WPIC}) if the complex sectional curvature is nonnegative for all these totally isotropic sections. 
This condition also originates in the work of Micallef-Moore \cite{MM}, where it arises naturally in the study of minimal surfaces.

By writing $v=a+ib$, with $a,b\in \R^n$ and expanding
$$g(v,v)=g(a,a)-g(b,b)+2ig(a,b),$$
we see that $v$ is isotropic if and only if $a$ and $b$ are orthogonal and of the same length. 
Given a totally isotropic section $\Si$ we can always find an orthonormal set 
$\{e_1,\ldots, e_4\}\subset \R^n$ such that $\Si$ is spanned by $e_1+ie_2$ and $e_3+ie_4$. 
To see this, take any orthogonal $v,w\in\Si$, each of length $\sqrt{2}$.
Each of them can be written $v=e_1+ie_2$ and $w=e_3+ie_4$ for orthonormal $e_1,e_2$ and orthonormal $e_3,e_4$, but we still need to show that $e_1,\ldots, e_4$ are orthonormal.
In addition to $v$ and $w$ being isotropic, any linear combination of $v$ and $w$, for example $v+w$, is also isotropic, and so
$$0=g(v+w,v+w)=g(v,v)+g(w,w)+2g(v,w)=2g(v,w)=2\langle v,\bar w \rangle,$$
and we see that $v$ is orthogonal to $\bar w$ as well as $w$.
Rephrased, we have $v$ orthogonal to both $e_3$ and $e_4$ as required.
In particular, we can only have totally isotropic sections $\Si$ if $n\geq 4$. 
We will see that positive isotropic curvature can be used to restrict the topology of a manifold
but the stronger topological conclusions we will be seeking require
strengthening of the curvature condition.

Equipped with the complexified viewpoint, we can return to give an alternative definition of 2-nonnegative curvature operator. Instead of considering isotropic vectors in $\C^n$ with respect to $g$, we exploit our inner product $\ci$ on $\Lambda^2\R^n$, which can be extended to $\Lambda^2\C^n$ by complex linearity, in order to consider isotropic elements of $\Lambda^2\C^n$. 
We can then define 
$\calr\in \mathcal{C}_B(\R^n)$ to have \textbf{2-nonnegative curvature operator} if
$$\calr(\om,\overline{\om})\geq 0\quad\text{ for all \textbf{isotropic} }\om\in \Lambda^2\C^n,$$
i.e. for all $\om$ with $\ci(\om,\om)=0$.
To see that this recovers the classical definition, we return to the earlier argument that $v\in\C^n$ is isotropic if and only if $v=a+ib$ for orthogonal $a,b\in \R^n$ of the same length. 
The analogous statement in our new setting is that $\om\in \Lambda^2\C^n$ is isotropic if and only if we can write $\om=\alpha+i\beta$, for orthogonal $\alpha,\beta\in \Lambda^2\R^n$ of the same length. 
Then 
$$\calr(\om,\overline{\om})=\calr(\alpha+i\beta,\alpha-i\beta)=\calr(\alpha,\alpha)+\calr(\beta,\beta).$$
Positivity of the right-hand side for every orthogonal $\alpha,\beta\in \Lambda^2\R^n$ of the same length is precisely the classical definition of 2-nonnegative curvature operator.

The discussion above leaves us in an ideal position to define the main curvature condition that concerns us in this survey, being a direct combination of the notions of 2-nonnegative curvature operator and of nonnegative complex sectional curvature.

\begin{defn}
We define 
$\calr\in \mathcal{C}_B(\R^n)$ to be \textbf{weakly PIC1} (or simply \textbf{WPIC1}) if 
$$\calr(\om,\overline{\om})\geq 0\text{ for all \textbf{simple and isotropic} }\om\in \Lambda^2\C^n.$$
\end{defn}
We refer to the complex sections $\Si$ corresponding to such $\om$ as PIC1 sections, or simply \emph{degenerate} sections.
An alternative way of describing such $\Si$ is that 
there exists $v\in \Si$ so that $\bar{v}$ is orthogonal to $\Si$,
or equivalently so that $g(v,w)=0$ for all $w\in \Si$. 
We can use this final viewpoint to see that WPIC1 implies WPIC.

The Ricci curvature of $\calr$ can be written as a sum of complex sectional curvatures corresponding to PIC1 sections. More precisely if $a_1\in \R^n$ is a unit vector for which we would like to express $\Ric(a_1,a_1)$, we can complete to an orthonormal basis
$a_1,\ldots,a_n$ and then consider the PIC1 sections represented by elements $\om_{jk}\in \Lambda^2\C^n$ of the form $a_1\wedge (a_j+i a_k)$ for $1<j<k\leq n$. 
Because 
$$\calr(\om_{jk},\overline{\om_{jk}})=\calr(a_1\wedge a_j,a_1\wedge a_j)+
\calr(a_1\wedge a_k,a_1\wedge a_k),$$
the sum of two sectional curvatures, we find that
$$\Ric(a_1,a_1)=\frac{1}{n-2}\sum_{1<j<k\leq n}\calr(\om_{jk},\overline{\om_{jk}}).$$
In particular, the WPIC1 condition always implies nonnegative Ricci curvature.

Our unorthodox definition of WPIC1 coincides with the various other earlier definitions.
Amongst these we remark that the original definition of Micallef-Moore \cite{MM} was that a Riemannian manifold is WPIC1 if and only if its Cartesian product with $\R$ (or $S^1$) is WPIC. 
Similarly a manifold is WPIC2 if and only if its Cartesian product with $\R^2$ is 
WPIC.\footnote{Brendle-Schoen \cite{BS} were the first to discuss products with $\R^2$ being WPIC; it was then later realised \cite{NW07} that this notion coincided with the earlier notion of nonnegative complex sectional curvature.}
These interpretations justify the notation PIC1 and PIC2.
See also \cite{BrendleBook, wilking2013}.

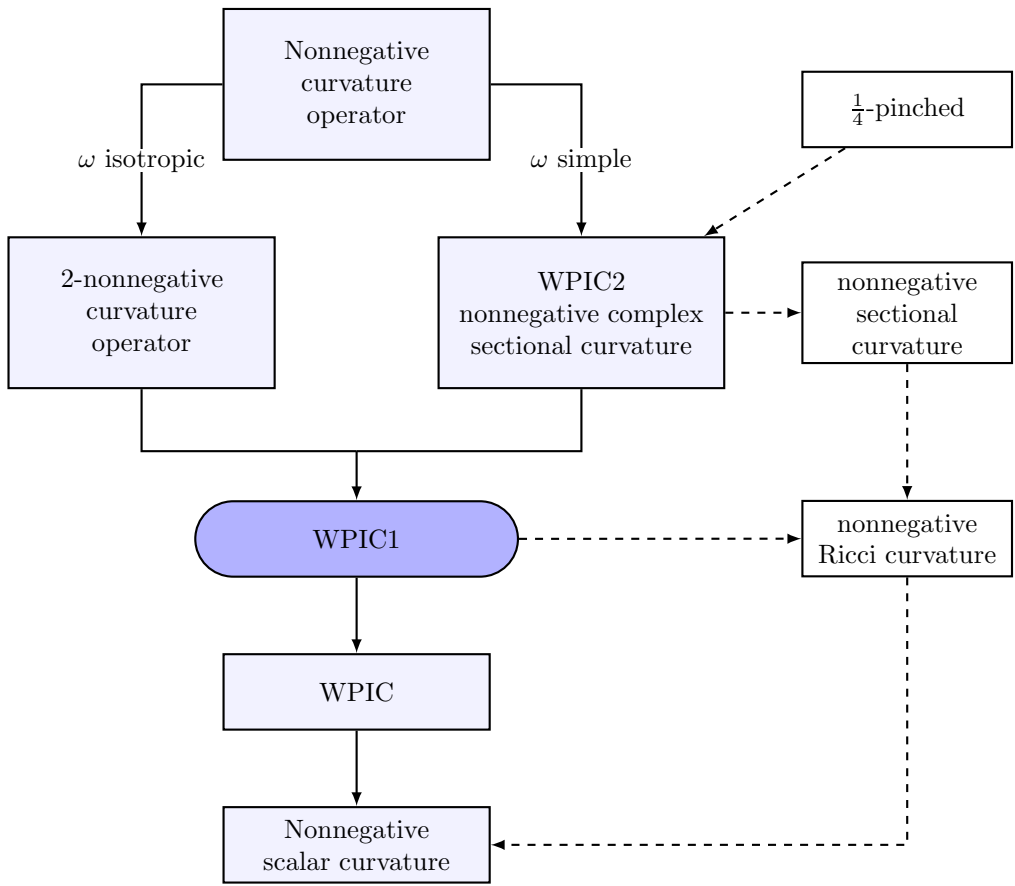
\begin{figure}
\begin{tikzpicture}[font=\small,thick]

\node[draw, fill=blue!5,
    rectangle, text width=2cm, text centered,
    minimum width=3.5cm,
    minimum height=2cm] (block4) {Nonnegative curvature operator};

\node[draw, fill=blue!5,
    rectangle, text width=2.5cm, text centered,
    below left=1cm and -0.7cm of block4,
    minimum width=3.5cm, 
    minimum height=2.0cm] (block5) {2-nonnegative curvature operator};
 
\node[draw, fill=blue!5,
    rectangle, text width=3.5cm, text centered,
    below right=1cm and -0.7cm of block4,
    minimum width=3.5cm,
    minimum height=2.0cm] (block6) {WPIC2\\ nonnegative complex sectional curvature};

\node[draw, text width=2.5cm, text centered,
    right=of block6,
    minimum width=2.5cm,
    minimum height=1cm] (block7) {nonnegative sectional curvature};

\node[draw, text width=2.5cm, text centered,
    above=1.5cm of block7,
    minimum width=2.5cm,
    minimum height=1cm] (block15) {$\frac14$-pinched};

\node[draw,
    rounded rectangle, fill=blue!30,
    below=4.5cm of block4,
    minimum width=4.5cm,
    minimum height=1cm,] (block11) { WPIC1};

\node[draw, text width=2.5cm, text centered,
    below=1.8cm of block7,
    minimum width=2.5cm,
    minimum height=1cm] (block8) {nonnegative Ricci curvature};

\node[coordinate,below=3.85cm of block4] (block12) {};

\node[draw,
    rectangle, fill=blue!5,
    below=of block11,
    minimum width=3.5cm,
    minimum height=1cm,] (block13) { WPIC};

\node[draw, text width=2.5cm, text centered,
    rectangle, fill=blue!5,
    below=of block13,
    minimum width=3.5cm,
    minimum height=1cm,] (block14) {Nonnegative scalar curvature};

\draw[-latex] (block4) -| (block5)
    node[pos=0.75,fill=white,inner sep=0]{$\om$ isotropic};
 
\draw[-latex] (block4) -| (block6)
    node[pos=0.75,fill=white,inner sep=0]{$\om$ simple};

\draw (block5) |- (block12);
\draw (block6) |- (block12);
\draw [-latex] (block12) -- (block11);
\draw [-latex] (block11) -- (block13);
\draw [-latex] (block13) -- (block14);

\draw [-latex, dashed] (block6) -- (block7);
\draw [-latex, dashed] (block11) -- (block8);
\draw [-latex, dashed] (block7) -- (block8);
\draw [-latex, dashed] (block8) |- (block14);
\draw [-latex, dashed] (block15) -- (block6);

\end{tikzpicture}
\caption{Map of curvature conditions}
\label{map_fig}
\end{figure}

The conditions we have seen can be summarised in the map of Figure \ref{map_fig}, along with some other important conditions and implications that will figure less prominently in this survey.

In general dimension, no two of the curvature conditions we have described on this map
coincide.
However in three dimensions 
all elements $\om\in\Lambda^2\R^n$ are simple, and so nonnegative curvature operator and nonnegative sectional curvature are equivalent notions. In turn these must then coincide with the intermediate notion WPIC2 on the map.
For the same reason, in three dimensions, 2-nonnegative curvature operator coincides with nonnegative Ricci curvature. Again, these will then coincide with the intermediate notion WPIC1 on the map.
The conditions WPIC1 and WPIC2 are thus distinct even in three dimensions because we can even find a homogeneous metric on $S^3$ that has positive Ricci curvature but does not have nonnegative sectional curvature. Indeed, 
we can consider a 
Berger sphere with $S^1$ fibres of length $2\pi t$.
When $t=1$ we have the standard  $S^3$.
More generally, the eigenvalues of the curvature operator are given by $t^2$, $t^2$ 
and $4-3t^2$.
The eigenvalues of the Ricci curvature are then given by $2t^2$, $4-2t^2$ and $4-2t^2$.
In particular, for $t^2> \frac{4}{3}$, the sectional curvature is no longer nonnegative,
but as long as $t^2< 2$, we still have positive Ricci curvature.

The equivalence of WPIC1 and nonnegative Ricci curvature in three dimensions is worth emphasising.
It is one of the morals of this survey that many statements that are true for nonnegative Ricci curvature in three dimensions will generalise to WPIC1 in higher dimensions.

\section{Preservation of positive curvature conditions under Ricci flow}

In order to understand the significance of the positive curvature conditions discussed in the previous section we would like to use Ricci flow. We can hope for this to work well
for curvature conditions that are preserved under the flow. 
If $M$ is a closed manifold and $g_0$ is a Riemannian metric on $M$ satisfying the curvature condition, then this amounts to asking that when we take a 
Ricci flow $g(t)$, $t\in [0,T)$, on $M$, i.e. a smooth family of Riemannian metrics satisfying 
$$\pl{g}{t}=-2\,\Ric_{g(t)},$$
then $g(t)$ also satisfies the curvature condition for $t\in [0,T)$.

Hamilton introduced a maximum principle technique known as the ODE-PDE theorem 
\cite{ham4PCO}, which demonstrates that this Ricci flow preservation on closed manifolds will always hold if the curvature condition is given by a closed, convex, $O(n)$-invariant 
set $F$ in $\mathcal{C}_B(\R^n)$ that is invariant under the so-called Hamilton ODE
\beq
\label{Ham_ODE}
\frac{d\calr}{dt}=Q(\calr),
\eeq
where $Q(\calr)$ is an algebraic curvature tensor that can be written in index notation as the quadratic expression
$$Q(\calr)_{ijkl}=R_{ijpq}R_{klpq}+2R_{ipkq}R_{jplq}-2R_{iplq}R_{jpkq}.$$
Many of the curvature conditions discussed so far are preserved under Ricci flow, 
including nonnegative curvature operator (Hamilton \cite{ham4PCO}), 
2-nonnegative curvature operator (Chen \cite{HChen}, Hamilton \cite{formations}), WPIC (H. Nguyen \cite{N10} and Brendle-Schoen \cite{BS}), and WPIC1 and WPIC2 (Brendle-Schoen \cite{BS}) in addition to 
nonnegative scalar curvature
(e.g. \cite[Corollary 3.2.3]{RFnotes}).
Aside from the easy case of nonnegative scalar curvature, 
which follows from an application of the classical parabolic maximum principle to the scalar curvature because of the evolution equation
\beq
\label{scalar_evol_eq}
\pl{\Scal}{t}=\lap\Scal+2|\Ric|^2
\geq \lap\Scal+\frac{2}{n}\Scal^2
\geq \lap\Scal,
\eeq
all of these are covered by the following result that is a 
rephrasing of a result of Wilking \cite{wilking2013}.

\begin{thm}
\label{wilking_interpretation}
Suppose $S\subset \Lambda^2\C^n$ is 
invariant under the action of $\mathrm{SO}(n,\C)$.
Then the convex cone
$$\cone(S):=\{\calr\in \CB\ :\ \calr(\om,\overline{\om})\geq 0\text{ for all }
\om\in S\}$$
is invariant under the Hamilton ODE. 
In fact, for all $h\in \R$, the set
$$G(S,h):=\{\calr\in \CB\ :\ \calr(\om,\overline{\om})\geq h\text{ for all }\om\in S\}$$
is invariant under the Hamilton ODE.
\end{thm}

To clarify, the action of $\mathrm{SO}(n,\C)$ on $\Lambda^2\C^n$ is determined by asking that if $A\in \mathrm{SO}(n,\C)$ and $v,w\in \C^n$, then 
$$A(v\wedge w)=Av \wedge Aw.$$

In most cases, the definitions of the preserved curvature conditions 
immediately give the corresponding sets $S$.
To summarise we have
\begin{enumerate}
\item
Nonnegative curvature operator: $S=\Lambda^2\C^n$,
\item
2-nonnegative curvature operator: $S=\{\om\in\Lambda^2\C^n\ :\ \om\text{ is isotropic}\}$,
\item
WPIC2: $S=\{\om\in\Lambda^2\C^n\ :\ \om\text{ is simple}\}$,
\item
WPIC1: $S=\{\om\in\Lambda^2\C^n\ :\ \om\text{ is simple and isotropic}\}$,
\item
WPIC: $S=\{\om\in\Lambda^2\C^n\ :\ \om\text{ is simple and every }u\in \Si(\om)
\text{ is isotropic}\}$,
\end{enumerate}
where if $\om=v\wedge w\neq 0$ for $v,w\in\C^n$, then $\Si(\om):=\mathrm{span}_\C\{v,w\}$.

In each of these cases, the definitions of $S$ only involve $g$, which remains invariant under the action of $\mathrm{SO}(n,\C)$ by definition, and so $S$ is automatically invariant as required.

In summary, we have a good selection of closed, convex, $O(n)$-invariant, ODE-invariant cones in $\CB$, and the ODE-PDE theorem to convert these to invariant sets for the curvature of Ricci flows on closed manifolds. It should, however, be stressed that on noncompact manifolds one can have a complete metric lying within an invariant cone, but a smooth Ricci flow evolution that no longer preserves this condition.

\section{The PIC1 sphere theorem}
\label{sphere_thm_sect}

Sphere theorems generally start with a closed Riemannian manifold that has some form of strictly positive curvature at each point, and deduce a global topological conclusion. 
A manifold that is PIC, i.e. that has positive complex sectional curvature on totally isotropic sections, has controlled topology, as demonstrated by Micallef-Moore
\cite{MM}. They showed that a closed, PIC, simply connected manifold is necessarily
homeomorphic to a sphere, generalising a number of classical sphere theorems. 
Without the simply connected hypothesis this fails: One can check that $S^1\times S^3$, with the standard metric, is PIC.

In this section we are looking for the weakest positive curvature condition that implies that a closed manifold is \emph{diffeomorphic} to a spherical space form, i.e. a quotient of a round sphere, without any hypothesis of simply connectedness.

Hamilton realised that Ricci flow could be used to address such problems.
He showed \cite{ham3D, ham4PCO} that closed manifolds of positive Ricci curvature in three dimensions, and positive curvature operator in four dimensions, are diffeomorphic to spherical space forms. In higher dimensions, Huisken \cite{huisken},
and independently Margerin \cite{margerin} and Nishikawa \cite{nishikawa}, arrived at the same conclusion with stronger curvature conditions,  and B\"ohm-Wilking \cite{BohmWilking2008} introduced new techniques for the analysis of solutions of the Hamilton ODE \eqref{Ham_ODE} in order to 
extend the theory to 2-positive curvature operator in arbitrary dimension. 
Brendle-Schoen \cite{BS} 
extended this theory to handle  PIC2 manifolds, 
i.e. manifolds of strictly positive complex sectional curvature,
and as a consequence established the  strict quarter pinching differentiable sphere theorem.

Given the map of curvature conditions in Figure \ref{map_fig}, all of the results just mentioned fit into the following extension.

\begin{thm}[Brendle \cite{brendlePIC1}, PIC1 sphere theorem]
\label{PIC1_sphere_theorem}
Suppose $(M,g_0)$ is a closed Riemannian manifold that is PIC1. 
Then $M$ is diffeomorphic to a spherical space form.
\end{thm}

Differentiable sphere theorems have already been written about extensively, e.g.
\cite{BrendleBook}, but we briefly survey the strategy.
The essential step is to construct a `pinching set', as introduced by Hamilton \cite{ham4PCO}. We give a definition that is almost equivalent to Hamilton's notion and the variants of B\"ohm-Wilking \cite{BohmWilking2008} and Brendle \cite{BrendleBook}.

\begin{defn}
A set $F\subset \CB $ is called a (generalised) pinching set if:
\begin{itemize}
\item
$F$ is closed, convex and $O(n)$-invariant;
\item
$F$ is invariant under the Hamilton ODE $\frac{d}{dt}\calr=Q(\calr)$;
\item
As $\la\downto 0$, 
\beq
\label{bd_cond}
\la F \to \R_{\geq 0}\ci.
\eeq

\end{itemize}
\end{defn}

The convergence in the final part is in the pointed Hausdorff sense. In other words, 
however large we take $r>0$, 
the intersection of $\la F$ with $B_r(0)\subset\CB$ converges in the Hausdorff sense to the intersection of $B_r(0)$ with the ray $\R_{\geq 0}\ci$ of curvature tensors of constant nonnegative sectional curvature.

The pinching set is specified so that if a Ricci flow has its curvature initially in $F$, then it remains so by the ODE-PDE theorem, and if the curvature blows up then the curvature is becoming like that of a sphere at each such point.
This motivates:

\begin{thm}[cf. Hamilton \cite{ham4PCO}]
\label{smooth_sphere_thm}
Let $(M^n,g_0)$ be a closed smooth Riemannian manifold, $n\geq 3$, with positive scalar curvature.
Suppose there exists a pinching set $F\subset \CB $ such that the curvature of $g_0$ lies in $F$ everywhere. 
Then $M$ is diffeomorphic to a spherical space form. 
Moreover, the maximal Ricci flow $g(t)$, $t\in [0,T)$, with $g(0)=g_0$ satisfies
$$\frac{1}{2(n-1)(T-t)}g(t)\to g_\infty,$$
smoothly as $t\upto T$, where $g_\infty$ has constant sectional curvature $1$.
\end{thm}

To clarify the term `maximal', 
Hamilton \cite{ham3D} showed that when the underlying manifold $M$ is closed, 
there exists a Ricci flow $g(t)$ on $M$, for $t\in [0,T)$, with either $T=\infty$, or with 
$$\sup_M |\calr|_{g(t)}\to\infty\text{ as }t\upto T.$$
This flow is unique in the sense that no other smooth Ricci flow with the same initial metric can deviate from $g(t)$ up to time $T$. It is called maximal because
we can then flow no further. Beware that on noncompact manifolds the curvature can blow up without the flow being maximal \cite{CabezasWilking2015, GT4}.

The flow of the theorem has 
positive scalar curvature. In particular, because $M$ is closed, the infimum of $\Scal(\calr_{g_0})$ is strictly positive.
The maximum principle applied to the parabolic equation \eqref{scalar_evol_eq} governing the scalar curvature then forces the scalar curvature to blow up, and thus the Ricci flow exists only on some finite time interval $[0,T)$.

Hamilton and Perelman \cite{RFnotes} tell us how as we approach a singular time $T<\infty$ in a Ricci flow on a closed manifold we can perform more and more extreme parabolic rescalings of the Ricci flow, 
cf. \cite[\S 1.2.3]{RFnotes},
and extract a limiting Ricci flow $g_\infty(t)$ on some possibly different underlying manifold $M_\infty$, for $t\in (-\infty,0]$, in an appropriate smooth Cheeger-Gromov-Hamilton sense. 
See \cite{RFnotes} for details.
The ODE-PDE theorem tells us that the curvature of $g(t)$ remains in $F$, and after rescaling to $g_\infty(t)$ the curvature must lie in $\lim_{\la\downto 0}\la F=\R_{\geq 0}\,\ci$ by definition of \emph{pinching set}.
In particular, at each point in spacetime $M_\infty\times (-\infty,0]$, all sectional curvatures must be the same, independent of which section we take. 
Schur's lemma then implies that at each time $t\in (-\infty,0]$, the manifold 
$(M_\infty,g_\infty(t))$ has constant positive sectional curvature, and thus must be a spherical space form. Such a manifold is necessarily compact, by Bonnet's theorem, 
and this forces $M$ to be diffeomorphic to $M_\infty$ by virtue of the properties of 
Cheeger-Gromov-Hamilton convergence. 

The sketched argument above implies the the first part of Theorem \ref{smooth_sphere_thm}, which is enough to prove sphere theorems.
The full claim concerning the smooth convergence of the Ricci flow as a time dependent tensor as $t\upto T$, and not just as $t_i\to\infty$ and requiring modification by diffeomorphisms, can be derived by combining with the work of Huisken
and Nishikawa \cite{huisken, nishikawa}.

\bigskip

Given what we have seen, the essential step in proving a sphere theorem such as
Theorem \ref{PIC1_sphere_theorem} is to find a suitable generalised pinching set $F$.
For the given metric $g_0$ on the given closed manifold $M$, we let $K$ be the set of all curvatures $\calr_{g_0}$, and then let $F\subset\CB$  be the  
smallest set containing $K$ that is closed, convex, $O(n)$-invariant and ODE invariant.
One then must prove that $F$ satisfies the blow-down condition \eqref{bd_cond}.

This was done by B\"ohm-Wilking \cite{BohmWilking2008} for 2-positive $(M,g_0)$
and Brendle-Schoen \cite{BS} gave the analogous result for PIC2, i.e. positive complex sectional curvature. 
Following Brendle \cite{BrendleBook}, to bootstrap the theory to PIC1 we need a  strengthening of the PIC2 result where one assumes some sort of almost-PIC2 condition, but this has to be done with some care because 
one can check
that complex projective space $\C P^n$ is on the boundary of the cone $\CPICt$ and being an Einstein manifold both the ODE and the PDE will evolve the curvature along the boundary of the cone $\CPICt$, 
by scaling, out to infinity without becoming round.
The resolution is to allow a notion of almost PIC2 that accommodates errors except where the curvature is large, when the curvature is expected to be well within the cone $\CPICt$. 
One way of asking that the curvature $\calr$ is a bit better than PIC2 for large 
$|\calr|$
is to ask for PIC2 pinching, i.e. that 
$$\calr-\ep\,\Scal(\calr)\,\ci\in \CPICt$$
for some small $\ep>0$.
Such $\calr$ make up a slightly thinner cone than $\CPICt$, more squeeezed towards the ray $\R_{\geq 0}\,\ci$.
An essentially equivalent way of squeezing the cone is to ask that 
$$\calr\in \ell_b(\CPICt)$$
for some small $b>0$, where $\ell_b$ is 
the endomorphism on $\CB$ introduced by 
B\"ohm-Wilking \cite{BohmWilking2008} defined by 
$$\textstyle \ell_b(\calr):=
\calr+b\cdot  \Ric(\calr)\owedge g+b^2\left(\frac{n-2}{n}\right) \Scal(\calr) \,  \ci.$$
Any compact $K$ in the interior of $\CPICt$ will also lie in the slightly smaller cone
$\ell_b(\CPICt)$ for some small $b>0$, and a result of 
B\"ohm-Wilking \cite[Proposition 3.2]{BohmWilking2008} tells us that $\ell_b(\CPICt)$ is also an ODE-invariant closed, convex, $O(n)$-invariant cone.

The Brendle-Schoen theory \cite{BS}, as in Brendle \cite{BrendleBook}, 
applies if we can find some (possibly small) $b>0$ and (possibly
large) $h\geq 0$ so that $F$ lies inside the set  of 
all $\calr\in\CB$ with 
\beq
\label{desired_F_location}
\calr+h\,\ci\in \ell_b(\CPICt).
\eeq
To recap, the $h\,\ci$ is giving us some freedom to lie outside the cone, whereas the $\ell_b$ is making the cone thinner; if $K$ lies in the interior of $\CPICt$, then we can take $h=0$.

To establish \eqref{desired_F_location} we need another ODE-invariant convex set
(not a cone) defined by 
$$G:=\{\calr\in \mathcal{C}_B(\R^n)\ :\ \calr(\om,\overline{\om})+|\ci(\om,\om)|\geq 0
\text{ for every simple }\om\in \Lambda^2\C^n\},$$
which coincides with the set $G$ considered by Brendle in different language \cite{BrendleBook}.
That the set $G$ is a preserved curvature condition
follows immediately from Theorem \ref{wilking_interpretation}  with
$$S=\{\om\in \Lambda^2\C^n\ :\ \om\text{ simple and }|\ci(\om,\om)|= 1\}$$
or simply
$$S=\{\om\in \Lambda^2\C^n\ :\ \om\text{ simple and }\ci(\om,\om)= 1\}$$
in which case
$$G=\{\calr\in \mathcal{C}_B(\R^n)\ :\ \calr(\om, \overline{\om})\geq -1
\text{ for every }\om\in S\}.$$
%
%
Meanwhile, from the definition of $G$ we see not only that 
$$\CPICt\subset G\subset \CPICo,$$
but also the scalings $\la G$ of $G$ interpolate between $\CPICt$ and $\CPICo$ 
as $\la$ varies in $(0,\infty)$ in the following sense.
The sets 
$$\la G=
\{\calr\in \mathcal{C}_B(\R^n)\ :\ \calr(\om,\overline{\om})+\la |\ci(\om,\om)|\geq 0
\text{ for every simple }\om\in \Lambda^2\C^n\}$$
are nested, getting larger as $\la$ gets larger. The intersection over all $\la>0$ (i.e. the limit as $\la\downto 0$ because we are working with closed sets) 
is $\CPICt$. In the other direction, 
$\la G$ exhausts the interior of $\CPICo$ as $\la\to \infty$.
In particular, the initial curvature $K$ lies in the interior of $\la G$ for some large $\la$.
By doing a once and for all scaling up of the initial $g_0$ (and hence scaling down of the set $K$), we may assume that $K$ lies in the interior of $G$.
Therefore  for small $b>0$, we can even say that $K$ lies in $\ell_b(G)$.

Next we record that
$$G\subset \{\calr\in \mathcal{C}_B(\R^n)\ :\ \calr+\ci\in \CPICt\},$$
because 
$$|\ci(\om,\om)|\leq \ci(\om,\overline{\om})$$
by Cauchy-Schwarz.


By definition of $\ell_b$ we have
$$\ell_b(\ci)=(1+2b(n-1)+b^2(n-2)(n-1))\,\ci=:h\,\ci,
$$
so for $\s\in G$ we have
$$\ell_b(\s) + h\, \ci = \ell_b(\s+\ci)
\in \ell_b(\CPICt).$$
%
%
Following B\"ohm-Wilking \cite{BohmWilking2008}, Brendle 
\cite[Proposition 11]{brendlePIC1}
tells us that $\ell_b(G)$ is ODE invariant. 
It is also a closed, convex, $O(n)$-invariant set. Thus 
$$F\subset \ell_b(G)\subset \{\calr\in\CB\ :\ \calr+h\,\ci\in \ell_b(\CPICt) \}$$
as required to invoke the PIC2 theory and conclude the PIC1 differentiable sphere theorem \ref{PIC1_sphere_theorem}.

\section{The Topology of open manifolds with positive curvature.}

Just as a sphere theorem tells us about the topology of a closed manifold with positive curvature of some form, we can expect similar curvature conditions to restrict the topology of \emph{open} manifolds, by which we mean complete noncompact manifolds. There are several important results of this type to survey, but this time, the refinement in this direction at the level of PIC1 can only be stated precisely as a conjecture.

Given an open (complete, noncompact, no boundary) surface $N$, if the negative part of the Gauss curvature is integrable then Cohn-Vossen \cite{CV} tells us that 
$$\chi(N)\geq \frac1{2\pi}\int_N K\,dV,$$
where $K$ is the Gauss curvature and $\chi(N)$ is the Euler characteristic.
In particular, if $K>0$ at each point then we deduce that $N$ is topologically the plane. 
In fact, Huber \cite{Huber1957} tells us that when the negative part of the Gauss curvature is integrable, the surface $N$ must be conformally a compact Riemann surface with finitely many punctures. Thus when $K>0$, the surface is biholomorphic to the plane.

The topological consequence of the two-dimensional theory generalises to higher dimensions via the theorem of Gromoll and Meyer.

\begin{thm}[{Gromoll-Meyer \cite{GM1969}}]
Any open manifold with positive sectional curvature is diffeomorphic to Euclidean space.
\end{thm}
In fact, the diffeomorphism will be given by the exponential map based at an appropriate point. This theory generalises in an elegant way to the case of nonnegative sectional curvature via the soul theorem (see Cheeger-Gromoll \cite{CG1972} and Perelman \cite{P1994}).

A natural question is whether a version of the Gromoll-Meyer theorem can hold for curvature positivity conditions in the map of Figure \ref{map_fig} other than positive sectional curvature. 
A famous result of Schoen-Yau addresses this in three dimensions using minimal surface techniques.
\begin{thm}[{Schoen-Yau \cite{SY1982}}]
Any open 3-manifold with positive Ricci curvature is diffeomorphic to $\R^3$. 
\end{thm}
We propose PIC1 as the correct analogue of positive Ricci curvature for the higher dimensional case.
\begin{conj}[{\cite{MFO2019}}]
Any open PIC1 manifold is diffeomorphic to Euclidean space.
\end{conj}
We have seen that PIC1 and positive Ricci curvature coincide in three dimensions, so our conjecture is consistent with the theorem of Schoen and Yau. Further evidence includes the work of He-Lee \cite{HL2021} that implies the conjecture in the case that the volume growth is maximal in the sense that for some $\ep>0$ and some (equivalently all) point $p$ in the manifold, we have $\VolB(p,r)\geq \ep r^n$, for all $r>0$, where $n$ is the dimension of the manifold.

Deviating from the Riemannian theme of this survey for a moment, there is another possible generalisation of the two-dimensional theory mentioned above to higher dimensions if we switch to the K\"ahler setting, as articulated in the following famous conjecture of Yau.
\begin{conj}[Yau]
Every open K\"ahler manifold with positive holomorphic bisectional curvature is biholomorphic to $\C^n$.
\end{conj}


\section{The topology of PIC1 pinched manifolds}
\label{pinch_sect}

In  Section \ref{sphere_thm_sect} we saw how sphere theorems turn curvature positivity hypotheses on closed manifolds into the deduction that the manifold is diffeomorphic to a space-form.
We now consider how more quantified versions of  curvature positivity on general complete manifolds can force the manifold to be closed.

The classical Bonnet-Myers theorem tells us that a complete manifold with 
quantifiably positive Ricci curvature in the sense that 
$\Ric\geq (n-1)r^{-2}>0$ for some $r>0$, i.e. all the eigenvalues of $\Ric$ are at least $(n-1)r^{-2}$, must have diameter no greater than $\pi r$. In particular it must be compact.

If we relax the curvature condition simply to $\Ric>0$ then the compactness conclusion fails because one can have a paraboloid, for example. Hamilton conjectured\footnote{Gerhard Huisken has informed us that Willmore was also discussing this conjecture in the early 1990s.
Hamilton wrote about the extrinsic version of the conjecture in 
\cite{Hamilton1994}. The first reference to the intrinsic case in the literature seems to be in \cite{ChenZhu2000}.}
\cite[Conjecture 3.39]{CLN09} that if one asks that the Ricci curvature is nonnegative, but in addition that at each point $p$ in the manifold the ratio of the largest and smallest eigenvalues of $\Ric$ is bounded by a constant that is independent of the point $p$,
then we can again draw topological conclusions. His conjecture has recently been proved by a combination of new work of Deruelle-Schulze-Simon \cite{DSS22} and of M.-C. Lee and the author \cite{LeeTopping2022}, incorporating earlier results and ideas of 
Lott \cite{Lott2019} and Chen-Zhu \cite{ChenZhu2000}.
See also the later paper of Huisken-K\"orber \cite{HK}.
\begin{thm}[{Hamilton's pinching conjecture,
cf. \cite[Conjecture 3.39]{CLN09}}]
\label{3D_thm}
Suppose $(M^3,g_0)$ is a complete (connected) three-dimensional Riemannian manifold with $\Ric\geq \ep\, \Scal\geq 0$ for some $\ep>0$. Then  $(M^3,g_0)$ is either flat or is compact.
\end{thm}
In higher dimensions there are many Ricci flat manifolds that are neither flat nor compact
(for example the Eguchi-Hanson space) so
we need to find a suitable replacement for Ricci pinching in higher dimensions.
Even more recent work together with M.-C. Lee \cite{LeeToppingPIC1pinch} illustrates how 
just as for sphere theorems, the PIC1 condition provides exactly what is required, as we now describe.

It is straightforward to compute that every $\calr\in\CB$ satisfies $\Scal(\calr)\leq c_n|\calr|$ for some constant $c_n$ depending only on the dimension $n$.
If we additionally assume that $\calr$ is WPIC1, then we know that 
$\Ric\geq 0$ and thus $\Scal\geq 0$. In fact, a short calculation \cite[Lemma A.2]{LeeToppingPIC1pinch} tells us that with this WPIC1 assumption, after possibly increasing $c_n$, we have 
$|\calr|\leq c_n\Scal(\calr)$.
In particular, every complex sectional curvature can be bounded from above in terms of the scalar curvature. 
The notion of PIC1 pinching is that every  complex sectional curvature 
corresponding to a PIC1 section
is \emph{comparable} to  the scalar curvature, i.e. it also has a lower bound in terms of the scalar curvature. This can be phrased as in the hypothesis \eqref{PIC1pinching} below.


\begin{thm}[{PIC1 pinching theorem, \cite{LeeToppingPIC1pinch}}]
\label{PIC1_pinching_thm}
Suppose $(M^n,g_0)$, $n\geq 3$, is a complete manifold of nonnegative complex sectional curvature 
that is PIC1 pinched in the sense that
\begin{equation}
\label{PIC1pinching}
\calr_{g_0}-\ep\,\Scal(\calr_{g_0})\, \ci\in \CPICo
\end{equation}
for some $\ep>0$.
Then $(M,g_0)$ is either flat or compact.
\end{thm}

The PIC1 pinching theorem generalises several earlier results in higher dimensions. 
Most recently, Brendle-Schoen \cite[Theorem 7.4]{BSsurvey} proved such a result in the positive scalar curvature case with the further hypotheses of PIC2 pinching (i.e. controlling all complex sectional curvatures to be comparable to the scalar curvature rather than just the PIC1 sections) and also with the assumption of uniformly bounded sectional curvature.
This in turn generalised earlier results of Ni-Wu \cite{NiWu2007} and Chen-Zhu \cite{ChenZhu2000}.

As an indicator of the heuristics behind Theorem \ref{PIC1_pinching_thm}, first observe that the nonnegativity of the Ricci curvature of $g_0$, as implied by $g_0$ being WPIC1, tries to make the manifold stabilise at infinity, as asserted by Gromov compactness.
Moreover, at least if $(M,g_0)$ has positive asymptotic volume ratio, 
we can blow down $(M,g_0)$ and extract a 
metric cone as a pointed Gromov-Hausdorff limit, as shown by 
Cheeger-Colding \cite{CC1}. At points near where the 
convergence is smooth,
the pinching condition will still hold on the limit. If $e_1$ is the radial direction of the cone, and we complete to a basis $e_1,\ldots,e_n$, then we have
PIC1 sections such as $e_1\wedge (e_2+ie_3)$ with zero complex sectional curvature, and the pinching then forces the scalar curvature and then the full curvature tensor to be zero at that point. If the whole cone away from the vertex arose as a smooth limit
then it would be a cone over a link of sectional curvature identically one.
In the case that the link is the unit sphere, this would imply volume growth identical to that of Euclidean space, also for the original manifold, and force the original manifold to be Euclidean space. For further discussion in the three-dimensional case, see \cite{Lott2019}.

This is very rough intuition that cannot be justified in all cases. However, it highlights a key problem of lack of regularity when we blow down, and hints that we might be able to solve it by employing the smoothing effects of Ricci flow. In fact, the most natural approach is to switch the blow-down intuition to the Ricci flow spacetime itself, as we will see.

\medskip

Given the Ricci flow technology that has been developed prior to \cite{LeeToppingPIC1pinch}, the key step in proving the PIC1 pinching theorem \ref{PIC1_pinching_thm} turns out to be a Ricci flow existence theorem. 
If we return to the discussion of sphere theorems in Section \ref{sphere_thm_sect}, we were always working on closed manifolds and so short-time Ricci flow existence theory followed trivially from the classical existence and uniqueness results of Hamilton \cite{ham3D}, simplified by De Turck \cite{deturck}.
On noncompact manifolds one does not expect to be able to find a complete Ricci flow 
starting with an arbitrary complete initial metric. However we make the following conjecture.
\begin{conj}[{\cite{MFO2019} and \cite[Conjecture 1.1]{MT2}}]
\label{WPIC1exist}
If $(M^n,g_0)$ is a complete WPIC1 manifold, $n\geq 3$, then there exists a complete WPIC1 Ricci flow $g(t)$ on $M$ for $t\in [0,T)$, some $T>0$, with $g(0)=g_0$.
\end{conj}
The three-dimensional version of this conjecture is simply saying that a complete initial metric with nonnegative Ricci curvature should induce a complete Ricci flow evolution. This case of the conjecture has been considered since the 1980s because its absence had led Shi's work \cite{shi_nonnegRic} on the classification of Ricci nonnegative three-manifolds to be stymied by the requirement of an additional hypothesis of bounded curvature. Although that problem was eventually solved by Liu \cite{gangliu} by switching to minimal surface techniques, Conjecture \ref{WPIC1exist} remains open even in three dimensions.
If one strengthens the hypothesis to WPIC2, then the desired solution is given by Cabezas-Rivas and Wilking 
\cite{CabezasWilking2015}.

In order to obtain the existence theory required to prove Theorem \ref{PIC1_pinching_thm}, it was necessary to exploit the pinching hypothesis. Together with M.-C. Lee, we proved:


\begin{thm}[{\cite[Theorem 1.3]{LeeToppingPIC1pinch}}]
\label{Thm:existence_mod}
For any $n\geq 4$, $\ep_0\in(0,\frac{1}{n(n-1)})$, there exist 
$a_0>0$ and $\ep_0'\in (0,\frac{1}{n(n-1)})$ such that 
the following holds. Suppose $(M^n,g_0)$ is a complete noncompact manifold 
such that
$$\calr_{g_0}-\ep_0\, \Scal(\calr_{g_0})\,\ci\in \CPICo$$
on $M$. Then there exists a smooth complete Ricci flow $g(t)$ on $M$
for $t\in [0,\infty)$, with $g(0)=g_0$, such that for all $t>0$, we have
\begin{enumerate}
\item[(a)] $\calr_{g(t)}-\ep_0' \Scal(\calr_{g(t)}) \,\ci\in \CPICo$;
\item[(b)] $|\calr|_{g(t)}\leq a_0 t^{-1}$.
\end{enumerate}
\end{thm}


The proof of this long-time existence theorem uses the Ricci flow theory that has been developed over recent years in order to handle noncompact or incomplete initial manifolds of potentially unbounded curvature, building on, in particular, the techniques from \cite{Hochard, ST1, ST2, hochard_thesis, cetraro}.
The intermediate step is to construct a local Ricci flow on a unit ball in 
$(M,g_0)$, for a time interval that only depends on $\ep_0$ and $n$ (and not otherwise on $g_0$) without boundary conditions, but with the substitute that the curvature 
$|\calr|_{g(t)}$ decays like $a_0/t$, which is a decay that is invariant under parabolic rescalings of the Ricci flow. The flow maintains a pinching condition with an error term.
The idea then is to apply the existence theorem not to $(M,g_0)$ itself but to blow-downs 
$(M,\la g_0)$ for smaller and smaller $\la>0$. The resulting local flows can be parabolically blown back up to match the initial metric $g_0$, rather than $\la g_0$, and we can take a limit of these blown up flows to give the global flow of the theorem. The process of rescaling destroys the error in the pinching estimate and  we recover exact PIC1 pinching as desired.
Further details can be found in \cite{LeeToppingPIC1pinch}.

\medskip

Equipped with the Ricci flow existence theorem \ref{Thm:existence_mod}, 
we can prove the PIC1 pinching theorem \ref{PIC1_pinching_thm}.
We outline the overall strategy 
in the case $n\geq 4$
and leave the details to \cite{LeeToppingPIC1pinch}, or to
\cite{LeeTopping2022} for the case $n=3$.

Suppose the assertion of Theorem \ref{PIC1_pinching_thm} is false.
Then we can find a complete noncompact $(M,g_0)$ of nonnegative complex sectional curvature, that is pinched in the sense of \eqref{PIC1pinching} but is not flat.
By lifting to the universal cover, we may assume that $M$ is simply connected.
By reducing $\ep>0$ if necessary (which just weakens the pinching hypothesis) we may assume that $\ep<\frac{1}{n(n-1)}$, so Theorem \ref{Thm:existence_mod}
applies to give a complete Ricci flow $g(t)$  for all time $t\in [0,\infty)$, starting at $g_0$,
that satisfies PIC1 pinching and $a_0/t$ curvature decay, i.e. conclusions (a) and (b) of that theorem. Moreover, the flow will inherit the additional property of being WPIC2 from $g_0$, and it is even possible to deduce that the sectional curvature becomes positive everywhere for positive times.

A theorem attributed to Gromoll-Meyer then tells us that 
by virtue of the curvature decay and the positive sectional curvature,
the injectivity radius can  be controlled by
$$\inj_{g(t)}\geq \de\sqrt{t},$$
for some $\de>0$. 
One consequence of this control is that we can control the volume of balls of radius $\sqrt{t}$ by
$$\VolB_{g(t)}(x,\sqrt{t})\geq \eta\, t^{n/2},$$
for some $\eta>0$.
By considering this information at larger and larger times, and transporting it back to time zero, it is possible to deduce that the asymptotic volume ratio of $g_0$,
and indeed of each $g(t)$,  must be the same value, and that value must be strictly positive. 

Equipped with this information we can then blow down $g(t)$ parabolically, i.e. for $\la_i\downto 0$ we consider 
the rescaled Ricci flows 
$$g_i(t):=\la_i\, g\left(\frac{t}{\la_i}\right).$$
These still satisfy the same $a_0/t$ curvature decay and injectivity radius lower bound, 
and according to Hamilton's compactness theorem \cite{RFnotes} this allows us to
pass to a subsequence and extract a limit Ricci flow $(M_\infty,g_\infty(t))$
of the same dimension, for $t>0$, also with positive asymptotic volume ratio.
Moreover, by an argument of Schulze-Simon \cite{SchulzeSimon2013} one can deduce that this Ricci flow is an expanding gradient soliton.
When coupled with the pinching, that is inherited in the limit, 
an argument of L. Ni \cite[Corollary 3.1]{Ni2005} forces $g_\infty(t)$ to be flat.
With care, this flatness can be transferred to the original metric $g_0$ giving a contradiction.

We have only attempted to give the outline of the proof of the PIC1 pinching theorem. Details can be found in \cite{LeeToppingPIC1pinch}.

\medskip

Although we are not giving here the details of the proof of the main Ricci flow existence theorem \ref{Thm:existence_mod} behind the PIC1 pinching theorem \ref{PIC1_pinching_thm}, 
a key point is to understand why the notion of PIC1  is so important at a technical level.
On the one hand, it is possible to prove that the property of being locally PIC1 pinched cannot deteriorate too rapidly within a Ricci flow that satisfies scale-invariant 
$|\calr|_{g(t)}\leq c_0 t^{-1}$ curvature decay. 
(See \cite[Lemma 3.6]{LeeToppingPIC1pinch}.)
The proof of that uses a localised form of Hamilton's so-called ODE-PDE theorem. 
On the other hand, we need to prove that in the presence of some form of local
PIC1 pinching, the noncompact Ricci flow regularises in the sense that it is 
\emph{forced} to satisfy local curvature decay of the form
\begin{equation}
\label{C_0_decay}
|\calr|_{g(t)}\leq C_0 t^{-1}
\eeq 
for some short time interval.

The proof of this latter regularisation result begins by imagining a contradicting sequence in which \eqref{C_0_decay} fails for larger and larger $C_0$ in a shorter and shorter time. Normally at this point one would make  appropriate parabolic rescalings and argue that we have compactness, and that a subsequence must 
converge to an ancient solution of Ricci flow that has absurd properties.
(For an argument of this form, see \cite[Lemma 2.1]{ST1}, incorporating ideas of Perelman \cite{Perelman2002}.)
For this particular application, in the light of work of Brendle-Huisken-Sinestrari \cite{BrendleHuiskenSinestrari2011} and Yokota \cite{TakumiYokota-CAG}, one might hope that an extracted limit would be a shrinking spherical space form, which would contradict the noncompactness of the original flow.

In practice, the argument is more delicate because when we blow up the sequence of Ricci flows, we do not have a suitable lower injectivity radius at our disposal, and so we cannot hope to have convergence of any subsequence. Instead, we must prove quantitative estimates on local regions of the blow-ups in order to establish that we get as close as we like to constant positive sectional curvature if we go far enough through the sequence of blown up Ricci flows. Eventually this quantitative roundness forces the manifold to close up on itself, contradicting the noncompactness of the flow as originally envisaged.
Details can be found in the proof of Lemma 3.6 in \cite{LeeToppingPIC1pinch}.


\section{PIC1 limit spaces}

We now turn to the topic of metric limit spaces, which at first glance may seem disconnected from the discussion so far, but turns out to have many common themes.
To simplify the discussion, we consider only  limit spaces 
$(X,d,p)$ that are complete
pointed Gromov-Hausdorff limits of complete pointed Riemannian manifolds $(M_i,g_i,p_i)$
of the same dimension $n\geq 3$, 
satisfying the properties
$$\left\{
\begin{aligned}
\calr_{g_i}+\alpha_0 \ci & \in\cone\\
\VolB_{g_i}(p_i,1) &\geq v_0
\end{aligned}
\right.$$
for constants $\al_0\in\R$ and $v_0>0$, where $\cone$ is a closed convex $O(n)$-invariant cone in $\CB$,
contained within the cone of all $\calr\in\CB$ with nonnegative Ricci curvature.
The volume lower bound is creating a \emph{noncollapsed} limit space, and in this case we may refer to it as being $n$-dimensional.
These limit spaces are often the additional spaces one must consider to make  classes of positively curved manifolds compact.

In the case that $\cone\subset \CB$ is the cone of elements having nonnegative sectional curvature, the limit spaces are instances of Alexandrov spaces. 
When $\cone$ is the cone of elements with nonnegative Ricci curvature, we obtain (noncollapsed) Ricci limit spaces, whose study was pioneered by Cheeger and Colding 
\cite{cheeger_book}.

Limit spaces can have singularities of various types. 
For example, if we take the Eguchi-Hansen space $(M^4,g)$, which is a complete Ricci-flat four-manifold that is asymptotic to the cone over $\R P^3$, then we can scale it down more and more to give a sequence of pointed manifolds $(M_i,g_i,p_i)=(M,\la_i g,p)$, $\la_i\downto 0$, as above that converge to 
its asymptotic cone at infinity, i.e. the cone over $\R P^3$. This Ricci limit space then has a singularity at the vertex where the tangent cone is not Euclidean space. In general, such singular points can be dense in the limit space. In this case, the singularity is particularly bad in that it is preventing the limit space from being a manifold: No neighbourhood of the vertex of the cone is homeomorphic to a Euclidean ball. 

A central question in this topic is to understand when the limit space \emph{does} have to be a manifold. It turns out that Ricci flow is extremely effective in answering this question owing to its regularisation properties, as foreseen by M. Simon \cite{simon2012}. Intuitively one would like to run the Ricci flow starting from each approximating space $(M_i,g_i)$, then pass to a limit of the flows to get a Ricci flow starting with the Ricci limit space in some sense. One can then hope to understand the structure of the limit space by comparing with the smooth Ricci flow.

In practice this is somewhat naive because we do not expect even short-time existence of the Ricci flow starting with manifolds as general as $(M_i,g_i)$, and even when we do get existence 
we do not generally have enough curvature decay of the solutions in order to get compactness of the flows and to construct a smooth limit solution.
However, this process can be carried out locally in three dimensions in order to prove the following, 
which is a little stronger than  a conjecture of Anderson, Cheeger, Colding and Tian in this dimension.
\begin{thm}[{\cite{ST2, Hochard, hochard_thesis}}]
In the case $n=3$, any (noncollapsed) Ricci limit space is locally bi-H\"older homeomorphic to a ball in $\R^3$.
\end{thm}

Of particular interest in this survey is the corresponding generalisation of this result to higher dimensions. The Eguchi-Hansen space demonstrates that \emph{Ricci limit space} is too weak a notion to make this work. Instead, Y. Lai \cite{yi_lai} showed that the proof can be carried over if we take PIC1 limit spaces, i.e. limit spaces $(X,d,p)$ as above with 
$\cone=\CPICo$, by virtue of the estimates of 
Bamler, Cabezas-Rivas and Wilking \cite{BCRW}; see also the  work of Hochard \cite{hochard_thesis}. The following result is a slight extension of this statement proved with McLeod \cite[Theorem 1.6]{MT2}.

\begin{thm}
\label{PIC1_limit_thm}
Any (noncollapsed, $n$-dimensional) PIC1 limit space is isometric to $(M,d,p)$, where $M$ is a smooth $n$-manifold, 
$p\in M$, and $d:M\times M\to [0,\infty)$ is a distance metric on $M$ such that 
for every smooth complete metric $g$ on $M$, the identity map $(M,d)\to(M,d_{g})$ is locally bi-H\"older.
\end{thm}

As above, the naive starting point to prove this theorem would be to try to construct a Ricci flow starting with the PIC1 limit space. However, 
as alluded to above, even in the case that the PIC1 limit space is a smooth manifold itself, we cannot expect such a Ricci flow to exist. 
Our solution, inspired by the \emph{partial Ricci flow} strategy of Hochard, was to construct a 
\emph{pyramid} Ricci flow starting with the PIC1 limit space in some sense, and which lives not on a cylindrical space-time region $M\times (0,T)$ but on some partial subset of such a space-time. 
More precisely, we construct a smooth manifold $M$ and an incomplete Ricci flow $g(t)$ on a subset 
$\Om\subset M\times (0,T)$, where the time slices 
$\Om_t:=\{x\in M\ :\ (x,t)\in \Om\}\subset M$
are nested:
$$\Om_s\supset \Om_t\quad\text{ for } 0<s<t<T$$
with $\Om_t$ exhausting $M$ as $t\downto 0$,
and where the Riemannian distance of $g(t)$ satisfies 
$$d_{g(t)}\to d\quad\text{ locally uniformly on }M\times M\text{ as }t\downto 0,$$
for some metric $d$ on $M$ for which we subsequently show that $(M,d)$ is isometric to the original PIC1 limit space.
Moreover, by invoking bi-H\"older estimates introduced in \cite[Lemma 3.1]{ST2}, 
we find that for every compact $K\subset M$, and $t\in (0,T)$ sufficiently small (in particular so that $K\subset \Om_t$), the identity map
$(K,d)\to (K,d_{g(t)})$ is bi-H\"older, where $d_{g(t)}$ is the Riemannian distance on $(\Om_t,g(t))$,
which is enough to prove Theorem \ref{PIC1_limit_thm}. A longer overview of techniques of this form can be found in 
\cite{cetraro}.

As in Section \ref{pinch_sect}, it is instructive to consider why the PIC1 condition is so effective in these results. At a technical level, we highlight two principles.
The first is that within a local Ricci flow satisfying scale-invariant regularisation
(see \eqref{reg_ests} below)
the least $\al_0\in \R$ for which 
$\calr_{g(t)}+\alpha_0 \ci $ is PIC1 cannot drop too fast as $t$ increases. 
The following result was proved for $n=3$ in \cite[Lemma 2.2]{ST1}, and generalised to higher dimensions in \cite[Proposition II.2.6]{hochard_thesis}, following \cite{BCRW}.
We state the variant from \cite[Lemma 3.2]{MT2}. See also \cite{LeeTam}.

\begin{lem}[{Propagation of local lower curvature bounds}]
\label{DB}
Let $n \geq 3$ and $c_0 , \al_0 > 0.$ Suppose that $( M , g(t) )$ is a smooth $n$-dimensional Ricci flow, 
defined for $t\in [0,T]$ and satisfying that for some point $x \in M$ and $\ep>0$ we have $B_{g(0)} (x , \ep) \subset \subset M.$ 
We further assume that
\beq
\label{reg_ests}
|\calr |_{g(t)} \leq \frac{c_0}{t} \qquad \text{and} \qquad \inj_{g(t)} \geq \sqrt{ \frac{t}{c_0} } 
\eeq
throughout $B_{g(0)} (x , \ep) \times (0, T]$ and that
$$\calr_{g(0)} + \al_0 \cal{I} \in \CPICo $$
throughout $B_{g(0)} (x,\ep)$.
Then there exist constants $S = S(n,c_0,\al_0,\ep) >0$ and 
$\al_1 = \al_1 ( n , c_0 , \al_0, \ep ) > 0$ such that
$$\calr_{g(t)} (x) + \al_1 \cal{I} \in \CPICo$$
for all times $0 \leq t \leq S\wedge T$.
\end{lem}

The second principle to highlight is that within a local Ricci flow for which 
$\calr_{g(t)}+\alpha_0 \ci $ is PIC1 throughout (for some new $\al_0$) we must have scale-invariant curvature decay, with constants depending on the initial volume. 
The following result was proved for $n=3$ in \cite[Lemma 2.1]{ST1},
\cite[Lemma 4.1]{ST2} and 
\cite[Lemma A.1]{MT1}, and generalised to higher dimensions in 
\cite[Lemma 3.4]{yi_lai}, following \cite{BCRW}.
We state the variant from \cite[Lemma 3.1]{MT2}.

\begin{lem}[Scale-invariant regularisation under local lower curvature bounds]
\label{loc_lemma_analogue}
Given  $n \in \N$ and $v_0 > 0$, there exists  
$C_0 \geq 1$ such that the following is true.
Let $\left( M , g(t) \right)$ be a smooth $n$-dimensional Ricci flow, $t \in [0,T]$, such that for some $p \in M$ and $\ep>0$ we have $B_{g(t)} (p,\ep) \subset \subset M$ for each $t \in [0,T],$
and so that for any $r \in (0,\ep]$ we have  
$\VolB_{g(0)} (p,r) \geq v_0 r^n $. 
Further assume that for some $\al_0 > 0$ and all $t \in [0,T]$
we have 
$$\calr_{g(t)} + \al_0 \mathcal{I} \in \CPICo \qquad \text{on} 
\qquad M\times [0,T].$$
Then there exists $S = S(n, v_0, \al_0, \ep) > 0$ such that for all 
$0 < t \leq S\wedge T$ we have 
$$|\calr|_{g(t)}(p) \leq \frac{C_0}{t} \qquad \text{and} 
\qquad \inj_{g(t)}(p) \geq \sqrt{ \frac{t}{C_0} }.$$
\end{lem}


Lemmata \ref{DB} and \ref{loc_lemma_analogue} are crucial tools in the building of both local and global Ricci flows in the setting of potentially unbounded curvature and/or rough initial data.


\medskip

As an addendum to this section we remark that techniques related to those in this section can be used to study the problem of when a smooth manifold that arises as a Gromov-Hausdorff limit 
of a 
sequence of smooth positively curved manifolds is itself positively curved despite the convergence not being $C^2$. Classical notions of positive curvature for which this problem has been considered include positive sectional curvature, for which the theory of Alexandrov spaces applies, and positive Ricci curvature, for which optimal transport theory can be invoked.
The WPIC1 elements of this story were investigated in \cite{LeeToppingIMRN2022}, where it was shown that when a sequence of WPIC1 Riemannian manifolds has a local Gromov-Hausdorff limit that is locally isometric to a smooth Riemannian manifold, then the limit is also WPIC1. See  \cite[Theorem 1.7]{LeeToppingIMRN2022} for a precise statement.

\bigskip

\noindent
\emph{Acknowledgements:} 
Thanks to ManChun Lee, Andrew McLeod and the Warwick geometric analysis reading group for useful conversations and comments.
This work was supported by EPSRC grant EP/T019824/1.
For the purpose of open access, the author has applied a Creative Commons Attribution (CC BY) licence to any author accepted manuscript version arising.

\vskip 0.2cm

\noindent
\url{https://homepages.warwick.ac.uk/~maseq/}

\noindent
{\sc Mathematics Institute, University of Warwick, Coventry,
CV4 7AL, UK.}

\end{document}